%% file: dme.tex
\documentclass{proc-l}

\usepackage{amssymb,amsfonts,amscd,xypic}
\usepackage{graphics}
\usepackage{epsfig}

\newtheorem{theorem}{Theorem}[section]
\newtheorem{lemma}[theorem]{Lemma}

\theoremstyle{definition}

\theoremstyle{remark}
\newtheorem{remark}[theorem]{Remark}

\copyrightinfo{2001}{American Mathematical Society}

\numberwithin{equation}{section}

\begin{document}


\title{On the Hartogs--Bochner phenomenon
for CR functions in $P_2(\mathbb{C})$}

\author{Roman Dwilewicz}

\address{Institute of Mathematics, 
Polish Academy of Sciences,
\'Sniadeckich 8, P.O.~Box 137, 00-950 Warsaw,
Poland}

\thanks{
Partially supported by a grant of the 
Polish Committee for Scientific Research 
KBN 2 PO3A 044 15 and by a grant from the French-Polish
program ``Polonium 1999''.}

\author{Jo\"el Merker}

\address{ Laboratoire d'Analyse, Topologie et Probabilit\'es,
Centre de Math\'ematiques et Informatique, UMR 6632,
39 rue Joliot Curie,
F-13453 Marseille Cedex 13, France. 
Fax: 00 33 (0)4 91 11 35 52}


\commby{Steve Bell}


\begin{abstract}
Let $M$ be a compact, connected, $\mathcal{C}^2$-smooth and {\it
globally minimal} hypersurface $M$ in $P_2(\mathbb{C})$ which divides
the projective space into two connected parts $U^{+}$ and $U^{-}$. We
prove that there exists a side, $U^-$ or $U^+$, such that every
continuous CR function on $M$ extends holomorphically to this
side. Our proof of this theorem is a simplification of a result
originally due to F. Sarkis.
\end{abstract}

\maketitle

\section{Introduction}

Historically, one of the amazing theorems in the theory of holomorphic
functions of several complex variables is the theorem of Hartogs of
1906 about extension of functions from a neighborhood of the boundary
of a domain to the inside of the domain. This theorem was generalized
by Bochner in 1943, namely that it is enough to consider
$\mathcal{C}^1$-smooth functions defined just on the boundary of the
domain that satisfy the tangential Cauchy-Riemann equations. Since
that time many versions and generalizations of the theorem appeared,
{\it see} \, for instance, Kohn-Rossi \cite{KR}, Ehrenpreis \cite{E},
Ivashkovich \cite{I}, Harvey \cite{H}, Harvey-Lawson \cite{HaLa},
Laurent-Thi\'ebaut \cite{La}, Dolbeault-Henkin \cite{DH}, Sarkis
\cite{S1}, \cite{S2}, and others - for a review see \cite{La}.

Recently, solving the boundary problem in the sense of Harvey and
Lawson for the graph of CR functions defined on boundaries of domains
in disc-convex K\"ahler manifolds, Sarkis \cite{S1} proved extension
of CR-meromorphic mappings with values in $P_2(\mathbb{C})$ (for
$P_n(\mathbb{C})$, $n\geq 3$, {\it see} \cite{HaLa}, {\it part
II}). As an application, he obtained Theorem 1.1 below. Our
contribution to this subject essentially lies in a simplification of
his proof, as will appear below.

Recall that a CR manifold $M$ (locally embeddable or not) is called
{\it globally minimal} if any two points of it can be joined by a
piece-wise smooth curve running in complex tangential directions
(\cite{Tr}). The result is the following:

\begin{theorem}
Let $M$ be a compact connected $\mathcal{C}^2$-smooth real
hypersurface in $P_2(\mathbb{C})$ that divides the projective space
into two open parts $U^-$ and $U^+$. If $M$ is globally minimal, then
\begin{enumerate}
\item[{\bf (1)}] There exists a side, $U^-$ \text{\rm or} $U^+$, to
which every continuous CR function on $M$ extends holomorphically.
\item[{\bf (2)}] All holomorphic functions on the other side of $M$
which are continuous up to $M$ are constant.
\end{enumerate}
\end{theorem}

\begin{remark}
Theorem 1.1 also holds true ({\it cf.} \cite{S1}, \cite{S2}) if,
instead of a globally minimal $M$ as above and instead of CR functions
on $M$, we consider an arbitrary $\mathcal{C}^2$-smooth hypersurface
as above and {\it holomorphic functions in a neighborhood
$\mathcal{V}$ of $M$ in $P_2(\mathbb{C})$}. Our proof is valid in
$P_n(\mathbb{C})$ for $n\geq 2$. Of course, the implication {\bf (1)}
$\Rightarrow$ {\bf (2)} is trivial.
\end{remark}

\medskip
\noindent
{\it Summary of the proof.} \ Let $f$ be a continuous CR function on
$M$. Using global minimality, we shall extend $f$ to a one-sided
neighborhood $\mathcal{V}^\pm(M)$ of $M$ in $P_2(\mathbb{C})$. By
deforming $M$ into $\mathcal{V}^\pm(M)$, we shall argue that we can
suppose that $M$ and $f$ are $\mathcal{C}^\infty$-smooth in the
assumptions of Theorem~1.1, and even that $f$ is holomorphic in a
neighborhood of $M$. As the Dolbeault cohomology group
$H^{0,1}(P_n(\mathbb{C}))$ vanishes for $n\geq 2$ ({\it see}
\cite{GH}, \cite{HeLe}), every $\mathcal{C}^\infty$ CR function $f$ on
$M$ can be decomposed as a jump $f = f^+ - f^-$ of some functions
$f^\pm$, holomorphic on $U^\pm$ and $\mathcal{C}^\infty$ on
$\overline{U^\pm}$. This decomposition property (which holds without
assuming global minimality) easily implies that points {\bf (1)} and
{\bf (2)} from Theorem 1.1 are equivalent in the $\mathcal{C}^\infty$
category ({\it see} Lemma~3.1). Afterwards, using a theorem of
Takeuchi \cite{T}, we may embed a strip neighborhood of $M$ in
$\mathbb{C}^N$ which then bounds a complex manifold to which $f^-$
{\it and} $f^+$ {\it both extend holomorphically} ({\it see}\cite{H}
and \cite{HaLa}, {\it part I}), and we easily deduce using the maximum
principle that either $f^-$ or $f^+$ is constant.

\medskip
\noindent
{\it Open question.} \ In Theorem 1.1, the question arises whether the
assumption that $M$ is globally minimal can be removed. This question
appears to have deep relations to some well-known conjectures of
foliation theory, especially the (non-)existence of Levi flat
hypersurfaces in projective spaces. Indeed, every compact CR manifold
$M\subset P_2(\mathbb{C})$ can be decomposed as a disjoint union of CR
orbits. Some of them are of dimension 3, and the others are of
dimension 2, i.e., Riemann surfaces, which cannot be closedly embedded
in $M$. The closure of each such Riemann surface defines a non-trivial
lamination of $P_2(\mathbb{C})$ and every other Riemann surface orbit,
having common points with the closure, is dense in this
lamination. {\it The open question is whether such laminations exist}
({\it cf.}~\cite{Gh}). For $n\geq 3$ (only), this question has been
fixed negatively for laminations arising from global foliations of
$P_n(\mathbb{C})$, for real analytic Levi-flat hypersurfaces ({\it
see} Cerveau \cite{C}, Lins Neto \cite{Li}) and for smooth Levi-flat
hypersurfaces ({\it see} Siu \cite{Si}). Thus, as a part of the
folklore, we conjecture that every smooth hypersurface $M\subset
P_n(\mathbb{C})$, $n\geq 2$, is globally minimal. \endremark

\medskip
\noindent
{\it Acknowledgement.} \
We are grateful to the referee for his clever support in fixing some
incorrections in the manuscript and for having pointed out to us substantial 
improvements in the proofs.

\section{Decomposition of CR functions}

As we are essentially interested in the geometric approach and not to
the best suited regularity assumptions, we shall restrict our
attention to the $\mathcal{C}^\infty$-smooth category. By a
deformation argument, we shall see in \S3 below that the case where
$M$ is $\mathcal{C}^2$ and the CR function $f$ on $M$ is only
continuous can be reduced to the case where both $M$ and $f$ are
$\mathcal{C}^\infty$. We recall that $H^{0,1}(P_2(\mathbb{C}))=0$ (for
$\mathcal{C}^\infty$-smooth forms). It is easy to deduce:

\begin{lemma}
Let $M$ be a $\mathcal{C}^\infty$ real hypersurface in
$P_2(\mathbb{C})$ dividing it into two open parts $U^-$ and $U^+$.
Every $\mathcal{C}^\infty$-smooth CR function $f$ on $M$ decomposes as
$f=f^+-f^-$, where $f^{\pm}\in \mathcal{O}(U^{\pm})\cap
\mathcal{C}^\infty(\overline{U^{\pm}})$.
\end{lemma}

\begin{proof}
We can extend $f$ to a $\mathcal{C}^\infty$-smooth function $F$ over
$P_2(\mathbb{C})$ in such a way that supp\,$F$ lies in an arbitrarily
small neighborhood of $M$ and $\overline \partial F|_M$ vanishes to
infinite order. We define
$\omega = \overline \partial F$ on $U^{+}$ and $\omega = 0$ 
on $U^{-}$.
The form $\omega$ is a $\mathcal{C}^\infty$-smooth $(0,1)$-form. As
$H^{0,1}(P_2(\mathbb{C}))=0$, we can solve the equation $\overline
\partial u = \omega$ with $u$ of class $\mathcal{C}^\infty$. So we
have
$$
\aligned 
& \overline \partial (F - u) = 0 \quad \text{on} \quad U^{+}
\quad \text{i.e., $F - u$ is holomorphic on $U^{+}$}, \\ & \overline
\partial u = 0 \qquad \qquad \text{on} \quad U^{-} \quad \text{i.e.,
$u$ is holomorphic on $U^{-}$,} 
\endaligned
$$
and furthermore
$$
F = (F-u) - (-u), \qquad F|_M = f.
$$
Obviously the components of the decomposition are of class
$\mathcal{C}^\infty$.
\end{proof}

\section{Embedding of a strip neighborhood of $M$ in $\mathbb{C}^N$}

Let $M$ be a $\mathcal{C}^\infty$-smooth hypersurface in
$P_2(\mathbb{C})$ bounding $U^-$ and $U^+$ as above. As a preliminary
to the proof of Theorem~1.1, we may establish an interesting
equivalence between {\bf (1)} and {\bf (2)} in the
$\mathcal{C}^\infty$ category.

\begin{lemma}
Let $M$ be a $\mathcal{C}^\infty$ real hypersurface in
$P_2(\mathbb{C})$ dividing it into two open parts $U^-$ and
$U^+$. Then the following properties are equivalent:
\begin{enumerate}
\item[{\bf (1')}] There exists a side, $U^-$ \text{\rm or} $U^+$, to
which every $\mathcal{C}^\infty$ CR function on $M$ extends
holomorphically.
\item[{\bf (2')}] All holomorphic functions on the other side of $M$
and $\mathcal{C}^\infty$ up to $M$ are constant.
\end{enumerate}
\end{lemma}

\begin{proof}
Without loss of generality, we can fix $U^+$ to be this side. Assume
therefore that every $\mathcal{C}^\infty$-smooth CR function on $M$
extends holomorphically to $U^+$. Let $g\in \mathcal{O}(U^-)\cap
\mathcal{C}^\infty (\overline{U^-})$. Then the trace $g\vert_M$
extends holomorphically to $U^+$. So $g$ extends holomorphically to
$P_2(\mathbb{C})$, whence it is constant; this is {\bf (2')}.

Conversely, suppose that $\mathcal{O}(U^-)\cap
\mathcal{C}^\infty(\overline{U^{-}})$ consists of constant functions
only. Let $f$ be a $\mathcal{C}^\infty$-smooth CR function on $M =
\partial U^-$. From Lemma~2.1 we obtain that $f = f^+ - f^-$, where
$f^+$ and $f^-$ belong to the evident spaces. By hypothesis, $f^-$ is
constant. Consequently, the CR function $f$ holomorphically extends to
$U^+$; this is {\bf (1')}.
\end{proof}

Let now $M$ be an orientable, $\mathcal{C}^2$-smooth, real hypersurface in
a complex manifold $X$. By $\mathcal{V}^\pm(M)$ we denote a {\it
one-sided} (global) {\it neighborhood} of $M$, {\it i.e.} a connected
open set which contains a (local) one-sided neighborhood of $M$ at
each point of $M$ and such that $\mathcal{V}^\pm(M) = \hbox{Int}
(\overline{\mathcal{V}^\pm(M)})$, so that $\mathcal{V}^\pm(M)$ contains a
neighborhood in $X$ of a point $p\in M$ whenever it contains the two
(local) one-sided neighborhoods at $p$. Such one-sided neighborhoods
are usually constructed by gluing analytic discs to $M$. Then the side
where the discs are lying may of course depend on the point $p\in M$
and it may vary effectively if $M$ is pseudoconvex somewhere and
pseudoconcave elsewhere. Without convexity assumption, we have\,:

\begin{theorem}
{\rm (\cite{Tr}, \cite{M}, \cite{J})}
Let $M$ be an oriented globally minimal
$\mathcal{C}^\infty$-smooth embedded hypersurface in a
complex manifold $X$. Then there exists a one-sided neighborhood
$\mathcal{V}^\pm(M)$, constructed by gluing small
analytic discs to $M$, to which all continuous CR functions on $M$ can be
holomorphically extended.
\end{theorem}

\begin{proof}[Proof of Theorem 1.1.]
Recall that {\bf (1)} implies {\bf (2)} trivially. We prove {\bf
(1)}. At first, we show that we can assume that $M$ and $f$ are
$\mathcal{C}^\infty$. So, let $M$ be a compact connected
$\mathcal{C}^2$ real hypersurface in $P_2(\mathbb{C})$ dividing
$P_2(\mathbb{C})$ into two connected open components $U^+$ and $U^-$
and let $f$ be a continuous CR function on $M$. By the theorem above,
$f$ extends holomorphically to a one-sided neighborhood
$\mathcal{V}^\pm(M)$ of $M$. Let $f$ still denote the resulting
extension. Let us smoothly deform $M$ into $\mathcal{V}^\pm(M)$. We
get a $\mathcal{C}^\infty$ hypersurface $M'\subset\subset
\mathcal{V}^\pm(M)$ arbitrarily close to $M$ and the trace
$f':=f\vert_{M'}$ of $f$ on $M'$ defines a $\mathcal{C}^\infty$ CR
function (even better, we could assume that $M'$ and $f'$ are real
analytic). Notice that as global minimality is a stable property
under sufficiently small perturbations, the deformed hypersurface $M'$
can be assumed to be globally minimal also. Again, $M'$ divides
$P_2(\mathbb{C})$ into two open parts ${U'}^-$ and ${U'}^+$. Suppose
that we can show that there exists a side, say ${U'}^+$, to which
every $\mathcal{C}^\infty$ CR function $f'$ on $M'$ extends
holomorphically. By inspecting the relative geometry of the three open
sets $U^+$, $\mathcal{V}^\pm(M)$ and ${U'}^+$, and using the fact that
$f$ is holomorphic in the whole of $\mathcal{V}^\pm(M)$, we deduce
that every continuous CR function $f$ on $M$ extends holomorphically
to $U^+$. In summary, it suffices for us to prove Theorem~1.1 in the
case where both $M$ and $f$ are $\mathcal{C}^\infty$-smooth.

Thus, let $M$ be a compact connected $\mathcal{C}^\infty$ real
hypersurface of $P_2(\mathbb{C})$ which divides it into two open parts
$U^-$ and $U^+$. Let $f$ be a $\mathcal{C}^\infty$ CR function on $M$
and let $f^+\in \mathcal{O}(U^+)\cap
\mathcal{C}^\infty(\overline{U^+})$, $f^-\in \mathcal{O}(U^-)\cap
\mathcal{C}^\infty(\overline{U^-})$ be given by Lemma~2.1 so that
$f=f^+-f^-$ on $M$. By virtue of Lemma~3.1, it suffices to show that
at least one of the two functions $f^+$ or $f^-$ is constant. We then
proceed by contradiction.

Assume that $f^+$ {\it and} $f^-$ are both nonconstant. By global
minimality, $f$, $f^+$ and $f^-$ can be holomorphically extended to a
one-sided neighborhood $\mathcal{V}^\pm(M)$ of $M$ constructed by
gluing discs. This open set is in fact a one-sided strip open set,
although it is not, in general, a tubular neighborhood of
$M$. However, we can easily deform $M$ into $\mathcal{V}^\pm(M)$ in a
$\mathcal{C}^\infty$ fashion, getting a $\mathcal{C}^\infty$
hypersurface $M_0\subset\subset \mathcal{V}^\pm(M)$. Furthermore, we
can include $M_0$ in a one-parameter family of manifolds $M_t$,
$|t|\leq \varepsilon$, satisfying
\begin{enumerate}
\item[{\bf 1.}] $M_t\cap M_{t'} = \emptyset$ for all $t\neq t'$,
\item[{\bf 2.}] $P_2(\mathbb{C})\backslash M_t = U_t^+\cup U_t^-$,
\item[{\bf 3.}] For $t<t'$, we have $M_t\subset\subset U_{t'}^-$ and
$U_t^- \subset\subset U_{t'}^-$. Also, $M_{t'}\subset\subset U_{t}^+$
and $U_{t'}^+ \subset\subset U_t^+$,
\item[{\bf 4.}] \ The strip $S:= \bigcup_{|t|< \varepsilon} M_t$ is an
open neighborhood of $M_0$ which is contained in $\mathcal{V}^\pm(M)$.
\end{enumerate}
The schematic picture of the geometric situation is as follows:

\begin{center}
\input one-sided.pstex_t
\end{center}

\noindent 
Notice that $M_0\subset\subset U_\varepsilon^-$. We shall apply the
following theorem due to Takeuchi over the domain $U_\varepsilon^-$.

\begin{theorem}
{\rm (\cite{T})} Let $U\subset P_2(\mathbb{C})$ be a domain. Then
either
\begin{enumerate}
\item[(i)] Holomorphic functions on $U$ are constant, \text{\rm or}
\item[(ii)] Holomorphic functions on $U$ separate points.
\end{enumerate}
\end{theorem}

\begin{lemma}
For every $\delta$ with $0 <\delta < \varepsilon$, the holomorphic
functions from $\mathcal{O}(U_\delta^-)$ separate points and give a
local coordinate system in a neighborhood of each point of
$U_{\delta}^-$.
\end{lemma}

\begin{proof}
According to Theorem~3.3 and since there exists by assumption the
nonconstant function $f^-|_{U_\varepsilon^-}\in
\mathcal{O}(U_\varepsilon^-)$, then $\mathcal{O}(U_\varepsilon^-)$
separate points. We can therefore find at least two functions $h_j\in
\mathcal{O}(U_\varepsilon^-)$, $j=1,2$, such that the zero locus of
the Jacobian of the mapping $h:=(h_1,h_2)$ is a proper, closed,
complex analytic subvariety $\Sigma$ of $U_\varepsilon^-$ of complex
dimension $\leq 1$. Of course, $h$ gives a local coordinate system at
each point of $U_\varepsilon^-$ not in $\Sigma$. Fix $\delta$ with $0<
\delta < \varepsilon$. Let now $p\in U_\delta^-\cap \Sigma$ be
arbitrary. Then after composing $h$ with an automorphism $A$ of
$P_2(\mathbb{C})$ arbitrarily close to the identity, we obtain that
$h\circ A$ gives a local coordinate system at $p$. Such an
automorphism $A$ moves $U_\varepsilon^-$ a little bit. As $A$ can be
chosen to be arbitrarily close to the identity, we can insure that the
domain of definition of $h\circ A$ still contains $U_\delta^-$. This
completes the proof.
\end{proof}

Let now $\eta$ with $0< \eta < \delta < \varepsilon$. The above lemma
implies that the manifold with boundary $M_\eta\cup
U_\eta^-\subset\subset U_\delta^-$ is embeddable into some complex
euclidean space $\mathbb{C}^N$, $N\in \mathbb{N}$, through an
embedding $\Phi: M_\eta\cup U_\eta^-\to \mathbb{C}^N$ whose
components are holomorphic functions defined in $U_{\delta}^-$. The
result of the embedding is a complex {\it manifold} \,
$\Phi(U_\eta^-)$ with boundary the three-dimensional maximally complex
CR manifold $\Phi(M_\eta)$. Notice that $\Phi$ also embeds into
$\mathbb{C}^N$ the closed strip $S':=\bigcup_{|t|\leq \eta}M_t$ which
provides a tubular neighborhood of $M_0$. By construction, for all
$\vert t\vert \leq \eta$, the $\Phi(M_t)$ are oriented maximally
complex three-dimensional CR manifolds in $\mathbb{C}^N$. The
orientation on $\Phi(M_t)$ is simply the push-forward of the
orientation on $M_t$. Further, $\Phi(M_t)$ bounds the complex manifold
$\Phi(U_t^-)$.

Finally, the holomorphic function $f^+$ induces a nonconstant
holomorphic function on the strip $\Phi(S')$, say $g^+ := f^+ \circ
\Phi^{-1}|_{\Phi(S')}$. Applying a generalization of Bochner's theorem
({\it see} \cite{H}), we deduce that the CR function
$g^+\vert_{\Phi(M_t)}$ extends holomorphically to the complex
submanifold $\Phi(U_t^-)$ for all $\vert t \vert \leq \eta$ (the
smoothness of $\Phi(U_t^-)$ is needed, otherwise we could only say
that the {\it graph} of $g^+\vert_{\Phi(M_t)}$ extends as a complex
analytic set by Harvey-Lawson's theorem). Thanks to the maximum
principle, we may now derive the desired contradiction easily.

Recall that by assumption, $f^+$ is nonconstant. Therefore, by the
maximum principle applied to $f^+$ in the projective space over $U^+$,
we have for all $t$ with $-\eta \leq t<0$ the strict inequality
$$
\sup_{z\in M_t} |f^+(z)| > \sup_{z\in M_0} |f^+(z)|.
$$
because $M_0\subset\subset U_t^+$ and $f^+$ is nonconstant. But on the
other hand, by the maximum principle applied in the complex euclidean
space $\mathbb{C}^N$ to the nonconstant holomorphic function
$g^+|_{\Phi(S')}$ which extends holomorphically to $\Phi(U_\eta^-)$,
we have for all $t$ with $-\eta \leq t<0$ the reverse strict
inequality
$$
\sup_{z\in M_t} |f^+(z)|= \sup_{z\in \Phi(M_t)} |g^+(z)| < \sup_{z\in
\Phi(M_0)} |g^+(z)|= \sup_{z\in M_0} |f^+(z)|.
$$ 
This gives the desired contradiction, which completes the proof of
Theorem 1.1.
\end{proof}

\bibliographystyle{amsplain}

\end{document}

%% file: one-sided.pstex_t
\begin{picture}(0,0)%
\epsfig{file=one-sided.pstex}%
\end{picture}%
\setlength{\unitlength}{3947sp}%
\begingroup\makeatletter\ifx\SetFigFont\undefined%
\gdef\SetFigFont#1#2#3#4#5{%
  \reset@font\fontsize{#1}{#2pt}%
  \fontfamily{#3}\fontseries{#4}\fontshape{#5}%
  \selectfont}%
\fi\endgroup%
\begin{picture}(4973,1819)(276,-3288)
\put(2751,-3257){\makebox(0,0)[lb]{\smash{\SetFigFont{8}{9.6}{\familydefault}{\mddefault}{\updefault}{\sc Figure 1.}}}}
\put(5249,-2365){\makebox(0,0)[lb]{\smash{\SetFigFont{8}{9.6}{\familydefault}{\mddefault}{\updefault}$M_0$ }}}
\put(5182,-2643){\makebox(0,0)[lb]{\smash{\SetFigFont{8}{9.6}{\familydefault}{\mddefault}{\updefault}$M_t, t<0$}}}
\put(5219,-2066){\makebox(0,0)[lb]{\smash{\SetFigFont{8}{9.6}{\familydefault}{\mddefault}{\updefault}$M_t, t>0$}}}
\put(5038,-1624){\makebox(0,0)[lb]{\smash{\SetFigFont{8}{9.6}{\familydefault}{\mddefault}{\updefault}$M$}}}
\put(2963,-2291){\makebox(0,0)[lb]{\smash{\SetFigFont{8}{9.6}{\familydefault}{\mddefault}{\updefault}$U^ -$}}}
\put(4634,-1803){\makebox(0,0)[lb]{\smash{\SetFigFont{8}{9.6}{\familydefault}{\mddefault}{\updefault}$U_0^+$}}}
\put(4435,-2038){\makebox(0,0)[lb]{\smash{\SetFigFont{8}{9.6}{\familydefault}{\mddefault}{\updefault}$U_0^ -$}}}
\put(1673,-2733){\makebox(0,0)[lb]{\smash{\SetFigFont{5}{6.0}{\familydefault}{\mddefault}{\updefault}$\mathcal{V}^\pm(M)$}}}
\put(2111,-2991){\makebox(0,0)[lb]{\smash{\SetFigFont{5}{6.0}{\familydefault}{\mddefault}{\updefault}$\mathcal{V}^\pm(M)$}}}
\put(276,-1744){\makebox(0,0)[lb]{\smash{\SetFigFont{8}{9.6}{\familydefault}{\mddefault}{\updefault}$P_2(\mathbb{C})$}}}
\put(283,-2340){\makebox(0,0)[lb]{\smash{\SetFigFont{8}{9.6}{\familydefault}{\mddefault}{\updefault}$U^+$}}}
\end{picture}

%% file: dme.bbl
\begin{thebibliography}{10}

\bibitem{C} D.~Cerveau,
\textit{Minimaux des feuilletages alg\'ebriques
de $P_n(\mathbb{C})$.}
Ann.~Inst.~Fourier (Grenoble) \textbf{43} (1993), 1535--1543.

\bibitem{DH} P.~Dolbeault and G.M.~Henkin,
\textit{Cha\^{\i}nes holomorphes de bord
donn\'e dans $P_n(\mathbb{C})$.}
Bull.~Soc. Math.~France \textbf{125} (1997), 383--446.

\bibitem{E} L.~Ehrenpreis,
\textit{A new proof and an extension of Hartogs' theorem.}
Bull.~Amer.~Math.~Soc. \textbf{67} (1961), 507--509.

\bibitem{F} B.~Fabre,
\textit{Sur l'intersection d'une surface de Riemann
avec des hypersurfaces alg\'e\-briques.}
C.R.~Acad.~Sci.~Paris S\'er. I Math\'ematiques 
\textbf{322} (1996), 371--376.

\bibitem{Gh} \'E.~Ghys,
\textit{Laminations par des surfaces de Riemann.}
Dynamique et g\'eom\'etrie complexe (Lyon, 1997), ix, xi, 49--95, 
Panor. Synth\`eses, 8, Soc. Math. France, Paris, 1999.

\bibitem{GH} P.~Griffiths and J.~Harris, \textit{Principles of Algebraic
Geometry.} Wiley Classics Library, John Wiley and Sons, 1994.

\bibitem{H} R.~Harvey,
\textit{Holomorphic chains and their boundaries.} 
Proceedings of Symposia in
Pure Mathematics, \textbf{XXX} (1977), vol1, 
AMS, Providence, RI, 309--382.

\bibitem{HaLa} R.~Harvey and B.~Lawson,
\textit{On boundaries of complex analytic varieties.} Ann. of Math., I:
\textbf{102} (1975), 233 - 290; II: \textbf{106} (1977), 213--238.

\bibitem{HeLe} G.M.~Henkin and J.~Leiterer,
\textit{Theory of functions on complex manifolds.}
Monographs in Mathematics, \textbf{79}, Birkh\"auser, 
Basel-Boston, Mass., 1984, 226~pp.

\bibitem{I} S.~Ivashkovich,
\textit{The Hartogs-type extension theorem
for meromorphic maps into compact K\"ahler manifolds.}
Invent.~Math. \textbf{109} (1992), 47--54.

\bibitem{J} B.~J\"oricke,
{\it Some remarks concerning holomorphically convex hulls and envelopes of
holomorphy.} Math.~Z. \textbf{218} (1995), 143--157.

\bibitem{KR} J.J.~Kohn and H.~Rossi,
\textit{On the extension of holomorphic functions
from the boundary of a complex manifold.}
Ann.~of Math. \textbf{81} (1965), 451--472.

\bibitem{La} C.~Laurent-Thi\'ebaut,
\textit{Ph\'enom\`ene de Hartogs-Bochner
dans les vari\'et\'es CR.}
Topics in Complex Analysis,
Banach Center Publications, Warszawa {\bf 31}
(1995), 233--247.

\bibitem{Li} A.~Lins Neto,
\textit{A note on projective Levi flats
and minimal sets of algebraic foliations.}
Ann. Inst. Fourier (Grenoble) \textbf{49} 
(1999), 1369--1385.

\bibitem{M} J.~Merker,
\textit{Global minimality of generic manifolds
and holomorphic extendibility of CR functions.}
Int. Math. Res. Not. (1994), no1, 329--343.

\bibitem{S1} F.~Sarkis,
\textit{CR meromorphic extension and the non
embedding of the Andreotti-Rossi CR structure in the projective space.}
Int. J. Math. \textbf{10} (1999), 897--915.

\bibitem{S2} F.~Sarkis,
\textit{Probl\`eme de Plateau complexe dans les vari\'et\'es kahl\'eriennes.}
Preprint, 1999.

\bibitem{Si} Y.T.~Siu,
\textit{Nonexistence of smooth Levi-flat hypersurfaces in 
complex projective spaces of dimension $\geq 3$.}
Ann. Math. (2) \textbf{151} (2000), no.3, 1217--1243.

\bibitem{T} A.~Takeuchi,
\textit{Domaines psudoconvexes infinis et
la m\'etrique riemannienne dans un espace projectif.}
J.~Math.~Soc.~Japan \textbf{16} (1964) 159--181.

\bibitem{Tr} J.-M.~Tr\'epreau,
\textit{Sur la propagation des singularit\'es
dans les vari\'et\'es CR.}
Bull. Soc. Math. France
\textbf{118} (1990), 403--450.

\end{thebibliography}
